\documentclass[a4paper,11pt]{article}

\usepackage{subfig}
\usepackage{pgfplots}
\usepackage{pgfplotstable}
\usepackage{mathtools}
\usepackage{multicol}
\usepackage{comment}
\usepackage{booktabs}
\pgfplotsset{compat=1.5}
\usepackage{amssymb}
\usepackage{url}
\usepackage{bm}

\usepackage{pdfsync}
\usepackage{float}
\usepackage{tabularx}
\usepackage{enumerate}
\usepackage{array}
\usepackage{xspace}

\usepackage{authblk}

\newcommand{\x}{\ensuremath{\boldsymbol{x}}}
\newcommand{\mupar}{\ensuremath{\boldsymbol{\mu}}}

\begin{document}

\title{Shape optimization through proper orthogonal decomposition with
  interpolation and dynamic mode decomposition enhanced by active subspaces}

\author[]{Marco~Tezzele\footnote{marco.tezzele@sissa.it}}
\author[]{Nicola~Demo\footnote{nicola.demo@sissa.it}}
\author[]{Gianluigi~Rozza\footnote{gianluigi.rozza@sissa.it}}

\affil[]{Mathematics Area, mathLab, SISSA, International School of Advanced Studies, via Bonomea 265, I-34136 Trieste, Italy}

\maketitle

\begin{abstract}
We propose a numerical pipeline for shape optimization in naval
engineering involving two different non-intrusive reduced order
method (ROM) techniques. Such methods are proper orthogonal
decomposition with interpolation (PODI) and dynamic mode decomposition
(DMD). The ROM proposed will be enhanced by active subspaces (AS) as a
pre-processing tool that reduce the parameter space dimension and
suggest better sampling of the input space.

We will focus on geometrical parameters describing the perturbation of
a reference bulbous bow through the free form deformation (FFD)
technique. The ROM are based on a finite volume method (FV) to simulate
the multi-phase incompressible flow around the deformed hulls.

In previous works we studied the reduction of the parameter space in
naval engineering through AS~\cite{tezzele2018dimension,demo2018isope}
focusing on different parts of the hull. PODI and DMD have been employed for
the study of fast and reliable shape optimization cycles on a bulbous
bow in~\cite{demo2018shape}.

The novelty of this work is the simultaneous reduction of both the
input parameter space and the output fields of interest. In particular
AS will be trained computing the total drag resistance of a hull
advancing in calm water and its gradients with respect to the input
parameters. DMD will improve the performance of each simulation of the
campaign using only few snapshots of the solution fields in order to
predict the regime state of the system. Finally PODI will interpolate
the coefficients of the POD decomposition of the output fields for a
fast approximation of all the fields at new untried parameters given
by the optimization algorithm. This will result in a non-intrusive
data-driven numerical optimization pipeline completely independent
with respect to the full order solver used and it can be easily
incorporated into existing numerical pipelines, from the reference CAD to the
optimal shape.
\end{abstract}



\section{Introduction}
\label{sec:intro}
In a shape optimization problem, we aim to find the shape --- among all the
admissible geometries --- that minimizes a certain objective function. In this
work we propose a novel approach to optimize the total resistance of a
ship hull advancing in calm water by deforming the original hull, a common
problem in the naval engineering field.

First we define the total resistance as the sum of the viscous and
lift forces acting on the hull.
Formally, our optimization problem can be expressed as
\begin{equation}
\min_{\forall \mupar \in D} f(\Omega, \mupar) = 
\min_{\forall \mupar \in D} \oint p \cos{\theta}\:\mathrm{d} \Omega_{\mupar}
+ \oint \tau_x\:\mathrm{d} \Omega_{\mupar} ,
\end{equation}
where the $D \subset \mathbb{R}^P$ is the parametric domain, $P$ the
number of parameters, $\Omega \in
\mathbb{R}^3$ is the reference hull domain, and $\Omega_{\mupar} =
\mathcal{M}(\Omega, \mupar)$ is the defomed hull. The morphing map
$\mathcal{M}(\cdot, \mupar): \mathbb{R}^3 \to \mathbb{R}^3$ we use in this work is
the free form deformation (FFD) and will be properly defined in Section~\ref{sec:ffd}. Examples of other
deformation techniques are radial basis functions (RBF)
interpolation~\cite{buhmann2003radial,morris2008cfd,manzoni2012model}, and
inverse distance weighting (IDW)
interpolation~\cite{shepard1968,forti2014efficient,BallarinDAmarioPerottoRozza2018}.  
The unknowns $p$ and $\tau_x$ denote respectively the
pressure and the $x$ component of the wall shear stress over the hull surface,
while $\theta$ is the angle between the flow direction and the
surface. The evaluation of the objective function requires a
numerical simulation of the flow around the ship, which has a high
computational cost.
The purpose of this work is beyond an analysis of the adopted full order
model for the fluid dynamics, we just provide a brief summary in order to
facilitate the understanding of the pipeline. We resolve the Reynolds-averaged
Navier Stokes (RANS) equations with the $k$--$\omega$ SST turbulence
model using a finite volume approach, a typical benchmark in industrial
hydrodynamics analysis. Such model deals very well with turbulent fluid, but at
high computational cost. Moreover, due to the complexity of the optimization
problem, we typically need many evaluations of the objective function to
converge to the optimal shape.

\begin{figure}[h!]
\centering
\includegraphics[trim=50 160 50 200, width=.7\textwidth]{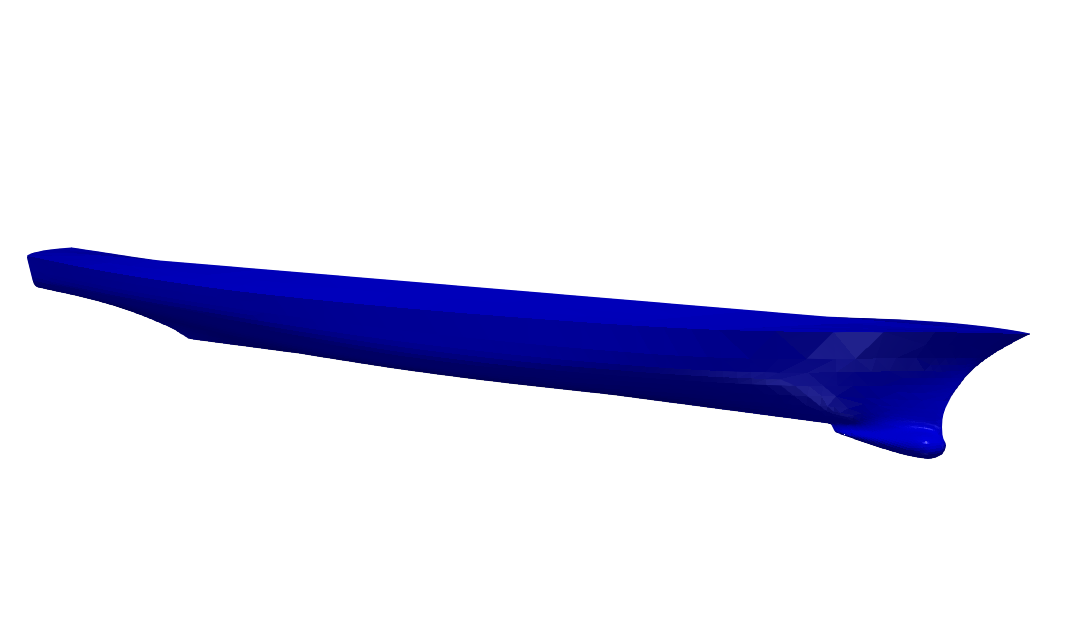}
\caption{Complete hull domain of the DTMB 5415.}
\label{fig:hull}
\end{figure}

For this work we choose to simulate the flow around the DTMB 5415 hull
due to the existence of a vast amount of literature and benchmark
tests. In Figure~\ref{fig:hull} the undeformed hull domain.

In order to reduce the computational cost, we introduce in the
optimization framework two reduced order modeling (ROM) techniques. These
techniqes are able to represent complex systems in a low dimensional space, reducing
the number of degrees of freedom used in the full order model discretization
and providing an efficient and reliable approximation of the solution. The ROM
methods initially collect a database of high-fidelity solutions --- the
solutions computed using the full-order model --- during the most
computationally expensive phase, also called {\it offline} phase. Then, the solutions
are combined to build the reduced space we query during the {\it online} phase
to obtain the new solution. In this work, we adopt the dynamic mode
decomposition (DMD) and proper orthogonal decomposition with interpolation
(PODI), two emerging data-driven techniques. 
PODI is used to approximate, given the high-fidelity solutions computed for
some defomed hulls, the solution for any new parametric point in the domain
$D$. DMD algorithm instead provides a simplification of the dynamics of complex
system: we use it in order to accelerate the single high-fidelity simulations
we need for PODI method, by storing few system outputs and exploiting them to
approximate the flow dynamics. For more details about equation-free ROM
methods, we suggest~\cite{tezzele2018ecmi}, while for a complete overview ---
including intrusive approaches --- we
cite~\cite{salmoiraghi2016advances,rozza2018advances,morhandbook2019}.

Moreover, additionally to these methods, we use the active subspace (AS)
property as preprocessing tool in order to be able to reduce the dimension
of the parameter space and obtain a better accuracy in ROM solution
approximation.

In this contribution, we focus on all the components of the
computational pipeline: in Section~\ref{sec:ffd} we provide a brief overview of
the FFD method, the Section~\ref{sec:dmd} illustrates the DMD
algorithm, in Section~\ref{sec:active} the AS property is explained,
while Section~\ref{sec:podi} describes the idea behind PODI technique.
Finally, Sections~\ref{sec:results}~and~\ref{sec:the_end} provide
respectively the numerical results collected during this work
and the final conclusions.

\section{The free form deformation technique}
\label{sec:ffd}

Free form deformationi (FFD) is a widespread deformation technique. Proposed
in~\cite{sederbergparry1986}, FFD was initially employed in computer graphics,
getting more popular both in academia and industry in the last decades. In this
section, we provide an overview of the method: for more details about FFD,
among all the works in literature, we recommend~\cite{rozza2013free,
demo2018shape, garotta2018quiet}.

\begin{figure}[h!]
\centering
\includegraphics[trim=0 0 0 0, width=.60\textwidth]{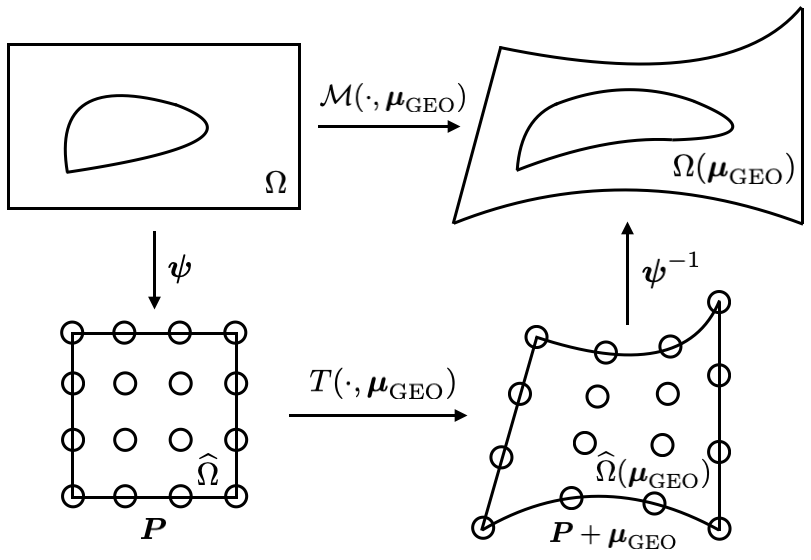}
\caption{Graphical representation of the FFD morphing map
  $\mathcal{M}$ as composition of the maps $\psi$, $T$, and
  $\psi^{-1}$. The displacements of the control points
  $\boldsymbol{P}$ define the morphing of the domain.}
\label{fig:ffd_scheme}
\end{figure}

\begin{figure}[h!]
\centering
\includegraphics[trim=0 70 140 0, width=.60\textwidth]{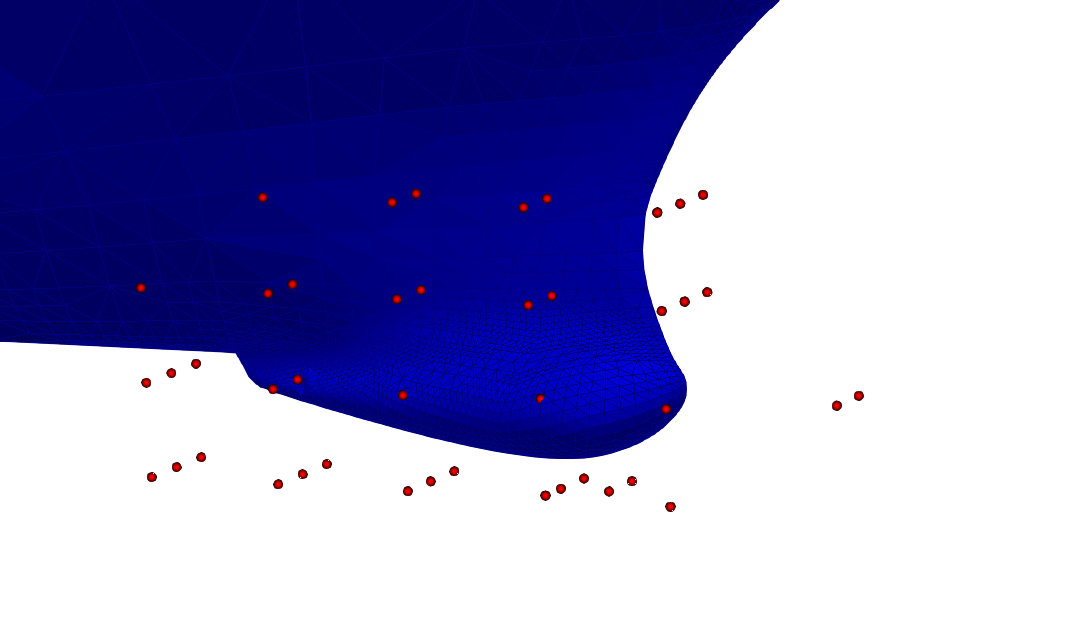}
\caption{Example of bulbous bow deformation using the FFD method. The red dots
are the FFD control points, already manipulated.}
\label{fig:ffd_lattice}
\end{figure}

The idea of FFD is very intuitive: the domain is deformed by manipulating a
lattice of points surrounding the object to morph. The displacements of these
control points are the input parameters $\mupar_{\text{GEO}}$. To achieve
this result $1$) the physical domain $\Omega$ is mapped to the
reference domain $\widehat{\Omega}$ using the function $\psi$, and a
lattice of control points $\boldsymbol{P}$ is constructed around the object to
deform, then $2$) through the map $T$ the reference domain 
is morphed using B-splines or Bernstein polynomials tensor product and finally
$3$) the deformed domain is remapped to the physical one by using $\psi^{-1}$.
In Figure~\ref{fig:ffd_scheme} is shown a sketch of the free form deformation
map as a composition of the three functions presented above.  Formally, we can
define the deformation map $\mathcal{M}$ as
\begin{equation}
\mathcal{M} (\x, \mupar_{\text{GEO}}) := (\boldsymbol{\psi}^{-1} \circ T
\circ \boldsymbol{\psi}) (\x, \mupar_{\text{GEO}})
\quad \forall \x \in \Omega.
\end{equation}
This technique allows to manipulate complex geometries and also computational
grids, since it is able to preserve derivatives continuity and perform global
deformation using only few parameters. Figure~\ref{fig:ffd_lattice} shows the
position of the lattice of control points around a bulbous bow, which
is the part of the hull we want to parametrize and morph. Regarding
the implementation, the results in this contribution are obtained
using PyGeM~\cite{pygem}, an open source Python package implementing
several deformation techniques. 

\section{Dynamic mode decomposition as accelerator of the single simulations}
\label{sec:dmd}

Dynamic mode decomposition is a data-driven modal decomposition
technique for analysing the dynamics of nonlinear
systems~\cite{schmid2010dynamic,schmid2011applications}. 
A comprehensive overview on DMD and its major variants
is in~\cite{kutz2016dynamic}. Other nonintrusive approaches with
randomized DMD can be found
in~\cite{bistrian2017randomized,bistrian2018efficiency}, while naval
engineering applications are in~\cite{demo2018shape,demo2018isope}.

Here we present a brief overview of the method and how we integrate it
in the computational pipeline we propose. Let us consider $m$
snapshots representing the state of the system for a given time
interval: $\{\boldsymbol{x}_i\}_{i=1}^{m} \in \mathbb{R}^n$. We seek a
linear operator $\mathbf{A}$ to approximate the nonlinear dynamics of
the state variable $\boldsymbol{x}$, that is $\boldsymbol{x}_{k+1}
= \mathbf{A} \boldsymbol{x}_k$. In order to find the DMD decomposition
we only need to approximate the eigenpairs of the operator
$\mathbf{A}$, without explicitly compute it. We proceed by dividing
the snapshots in two matrices $\mathbf{X}$ and $\mathbf{Y}$ as in the
following:
\begin{equation*}
\label{eq:matarranged}
\mathbf{X} =
 \begin{bmatrix}
  x_1^1   & x_2^1  & \cdots & x_{m-1}^1 \\
  x_1^2   & x_2^2  & \cdots & x_{m-1}^2 \\
  \vdots  & \vdots & \ddots & \vdots    \\
  x_1^n   & x_2^n  & \cdots & x_{m-1}^n 
 \end{bmatrix},\quad\quad
 \mathbf{Y} =
 \begin{bmatrix}
  x_2^1   & x_3^1  & \cdots & x_m^1  \\
  x_2^2   & x_3^2  & \cdots & x_m^2  \\
  \vdots  & \vdots & \ddots & \vdots \\
  x_2^n   & x_3^n  & \cdots & x_m^n 
 \end{bmatrix}.
\end{equation*}
With this representation we seek $\mathbf{A}$ such that $\mathbf{Y}
\approx \mathbf{A} \mathbf{X}$. Using the Moore-Penrose pseudo-inverse
operator, denoted by $^\dagger$, we express the best-fit matrix as
$\mathbf{A} = \mathbf{Y} \mathbf{X}^\dagger$.
We can compute the POD modes of the matrix $\mathbf{X}$ and project
the data onto the subspace defined by them. We use the truncated
singular value decomposition obtaining $\mathbf{X} \approx \mathbf{U}_r
\boldsymbol{\Sigma}_r \mathbf{V}^*_r$, where the unitary matrix
$\mathbf{U}_r$ contains the first $r$ modes. With these modes we can
compute the reduced operator  $\mathbf{\tilde{A}} \in
\mathbb{C}^{r\times r}$ as
$\mathbf{\tilde{A}} = \mathbf{U}_r^* \mathbf{A} \mathbf{U}_r =
\mathbf{U}_r^* \mathbf{Y} \mathbf{X}^\dagger \mathbf{U}_r =
\mathbf{U}_r^* \mathbf{Y} \mathbf{V}_r \boldsymbol{\Sigma}_r^{-1} \mathbf{U}_r^* \mathbf{U}_r =
\mathbf{U}_r^* \mathbf{Y} \mathbf{V}_r \boldsymbol{\Sigma}_r^{-1}$,
without the explicit computation of the full operator $\mathbf{A}$.
The reduced operator describe the evolution of the low-rank
approximated state $\boldsymbol{\tilde{x}}_k \in \mathbb{R}^r$ as
$\boldsymbol{\tilde{x}}_{k+1} = \mathbf{\tilde{A}}
\boldsymbol{\tilde{x}}_k$. We can then recover the high-dimensional
state $\boldsymbol{x}_k$ using the POD modes already computed:
$\boldsymbol{x}_k = \mathbf{U}_r \boldsymbol{\tilde{x}}_k$.

Using the eigendecomposition of the matrix $\mathbf{\tilde{A}}$, that
is $\mathbf{\tilde{A}} \mathbf{W} = \mathbf{W} \boldsymbol{\Lambda}$,
we are able to compute the eigenpairs of the full operator
$\mathbf{A}$. In particular the eigenvalues in $\boldsymbol{\Lambda}$
correspond to the nonzero eigenvalues of $\mathbf{A}$, while the
eigenvectors $\boldsymbol{\Phi}$ of $\mathbf{A}$ can be computed in two ways: by
projecting the low-rank approximation $\mathbf{W}$ on the
high-dimensional space $\boldsymbol{\Phi} = \mathbf{U}_r \mathbf{W}$,
or by computing them exactly with $\boldsymbol{\Phi} =
\mathbf{Y}\mathbf{V}_r \boldsymbol{\Sigma}_r^{-1} \mathbf{W}$. 

\begin{figure}[!htbp]
\centering\includegraphics[width=.7\textwidth]{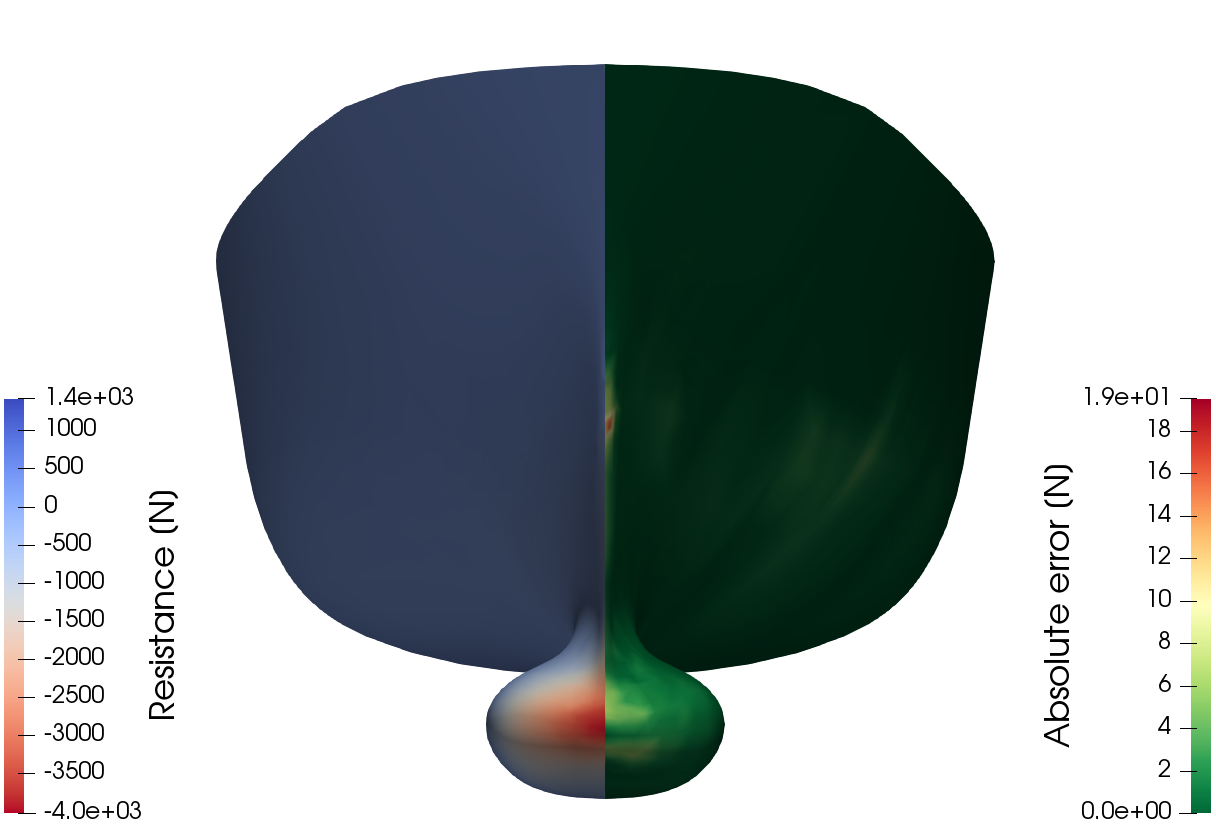}
\caption{Example of total drag reconstructed using the DMD algorithm. On the
left the reconstructed field, while on the right the absolute error with respect to
the high-fidelity solution.}
\label{fig:dmd}
\end{figure}
The actual implementation of the DMD algorithm we used and many
different variants from multiresolution
DMD~\cite{kutz2016multiresolution}, to DMD with
control~\cite{proctor2016dynamic}, and higher order
DMD~\cite{le2017higher}, can be found in the open source Python
package PyDMD~\cite{demo2018pydmd}.

In this work the DMD is used to accelerate the computation of the
total drag resistance for a given deformed hull. It uses only few
snapshots of the high-fidelity simulation, equispaced in time, to
predict the evolution of the target output. In particular we are
interested in the value of the total drag at regime. Figure~\ref{fig:dmd}
reports the forces field approximated using DMD and the absolute error
with respect to the high-fidelity solution for a particular
geometrical configuration.

\section{How to reduce the parameter space dimension with active subspaces}
\label{sec:active}

The active subspaces (AS) property has been formalized by Constantine
in~\cite{constantine2015active,constantine2014active}. It is a
property of a scalar function $f : \mathbb{R}^N \to \mathbb{R}$ and a
probability density function $\rho : \mathbb{R}^N \to \mathbb{R}^+$,
where $N$ is the number of the input parameters. Taking linear
combinations of the original parameters we can approximate $f$ using
these new parameters, thus reducing the parameter space dimension.
The output of interest $f(\mupar_{\text{GEO}})$, in our
case the total drag of the hull advancing in calm water, depends on
the geometrical parameters introduced in Section~\ref{sec:ffd}, while
$\rho$ describe the uncertainty in the model inputs, i.e. how we
sample the parameter space. For sake of clarity we will drop the pedix
and from now on $f(\mupar) := f(\mupar_{\text{GEO}})$. The general idea is to rotate the input
domain, after a proper rescale, in order to unveil a low dimensional parametrization of $f$,
which means to find proper directions in the input space where $f$
varies the most on average. We do so by checking the gradients of the
output of interest with respect to the parameters.

To proper exploit the AS property we introduce some
hypotheses: $f$ has to be continuous and differentiable with
square-integrable partial derivatives in the support of $\rho$. Then we
introduce the uncentered covariance matrix $\mathbf{C}$ of the
gradients of the target function, which is the matrix constructed with
the average products of partial derivatives of the map $f$ as follows

\begin{equation}
\label{eq:uncentered}
\mathbf{C} = \mathbb{E}\, [\nabla_{\mupar} f \, \nabla_{\mupar} f
^T] =\int (\nabla_{\mupar} f) ( \nabla_{\mupar} f )^T
\rho \, d \mupar,
\end{equation}
where with $\mathbb{E}$ we identify the expected value, and $\nabla_{\mupar} f =
\nabla f(\mupar) = \left [ \frac{\partial f}{\partial \mu_1}, \dots,
  \frac{\partial f}{\partial \mu_p} \right ]^T$ is the column vector
of partial derivatives of $f$. Since $\mathbf{C}$ is symmetric we can
express it with its real eigenvalue decomposition $\mathbf{C} =
\mathbf{W} \mathbf{\Lambda} \mathbf{W}^T$, where $\mathbf{W}$ is the
eigenvectors matrix, and $\mathbf{\Lambda}$ the diagonal matrix with
the eigenvalues in descending order. It can be proven that the
eigenvalues express the amount of variance of the gradient along the
corresponding eigenvector direction. This means that taking the first
$M$ most energetic eigenvalues and the corresponding eigenvectors, we
can approximate the target function with a reduce number of input
parameters. So the eigenpairs of $\mathbf{C}$ define the active
subspaces of the pair $(f, \rho)$. We proceed by partitioning
$\mathbf{W}$ and $\mathbf{\Lambda}$ as follows
\begin{equation}
\mathbf{\Lambda} =   \begin{bmatrix} \mathbf{\Lambda}_1 & \\
                                     &
                                     \mathbf{\Lambda}_2\end{bmatrix},
\qquad
\mathbf{W} = \left [ \mathbf{W}_1 \quad \mathbf{W}_2 \right ],
\end{equation}
where the pedix $1$ means the first $M$ eigenvalues and eigenvectors
respectively. Now we can use $\mathbf{W}_1$ to project the original
parameters to the active subspace, that is the span of the first $M$
eigenvectors. This means to align the input parameter space to
$\mathbf{W}_1$ and retain only the directions where $f$ varies the
most on average. We call active variable $\mupar_M$ the range of
$\mathbf{W}_1^T$, that is $\mupar_M = \mathbf{W}_1^T\mupar \in
\mathbb{R}^M$.
We can thus introduce a lower-dimension approximation
$g : \mathbb{R}^M \rightarrow \mathbb{R}$ of the quantity of interest $f$,
which is a function of $\mupar_M$ as follows
\begin{equation}
  f (\mupar) \approx g (\mathbf{W}_1^T \mupar) = g(\mupar_M).
\end{equation}

Active subspaces have been proven useful in naval applications
in~\cite{tezzele2018dimension,tezzele2018model,demo2018isope}, but
also coupled with POD-Galerkin model order
reduction~\cite{tezzele2018combined}. A gradient-free algorithm for
the discovery of active subspaces has been proposed
in~\cite{coleman2019gradient}, while an AS variant using average
gradients in~\cite{lee2019modified}.

We are going to find the active subspace for the total drag resistance
of the deformed hulls obtained by the FFD method and the application
of the DMD algorithm. Then we are going to exploit this active
subspace to perform a better sampling of the parameter space and thus
enhancing the construction of the reduced order model.

\section{Proper orthogonal decomposition with interpolation}
\label{sec:podi}
Reduced order modeling (ROM) is a popular technique to reduce the computational
cost of numerical simulations. Among all the available methods to achieve this
reduction, we focus in this contribution to the reduced basis method using the
proper orthogonal decomposition (POD) algorithm for the basis identification.
This method allows to reduce the number of degrees of freedom of a parametric
system by collecting the {\it snapshots} --- the full order system outputs ---
for several different configurations and combining them in an efficient way
for a real-time approximation of new solutions (for any new configuration).  In
the POD reduction framework, we can discern two main techniques:
POD-Galerkin, which requires all the details of the full order system
to generate a consistent low-dimensional representation of the physical problem, and
POD with interpolation (PODI), which instead requires only the snapshots. Due to these
requirements, the PODI method is particularly suited for industrial problem,
since it is able to been coupled to all the numerical solvers, even
commercial ones. In this contribution, we adopt PODI method. For more information
about POD-Galerkin, we
suggest~\cite{stabile2018finite,StabileHijaziLorenziMolaRozza2017,KaratzasBallarinRozza2019,georgaka2018parametric},
while for other examples of PODI applications we
recommend~\cite{garotta2018quiet,demo2018ezyrb, salmoiraghi2018}.

To calculate the POD modes we use the singular value decomposition (SVD) applied to the
snapshots matrix $\mathbf{X}$ such that $\mathbf{X} = \mathbf{U}
\mathbf{\Sigma} \mathbf{V}^*$. The columns of the unitary matrix $\mathbf{U}$
are the POD modes end the corresponding singular values, the elements in the
diagonal matrix $\mathbf{\Sigma}$ in decresing order, indicate the
energy associated to each mode.  Hence it is possible to select the first modes --- the most
energetic --- to span the reduced space and project onto it the high-fidelity
snapshots. In matricial form, we have:
\begin{equation}
\mathbf{X}^{\text{POD}} = \mathbf{U}_N^T \mathbf{X} ,
\end{equation}
where $\mathbf{U}_N$ is the matrix containing the first $N$ modes, and
$\mathbf{X}^{\text{POD}}$ is the matrix whose columns $\mathbf{x}^{\text{POD}}_i$ are the reduced
snapshots. We note that $\mathbf{x}^{\text{POD}}_i \in V^N$
and $\mathbf{x}_i \in V^\mathcal{N}$ where $\mathcal{N}$ refers to
the number of degrees of freedom of the full-order system.  Finally, due to the
reduced dimension, we are able to interpolate the reduced snapshots in order to
approximate the solution manifold. The new interpolated reduced snapshots are
then mapped back to the high-dimensional space for a real-time evaluation of
the solution.
To perform the non-intrusive model order reduction, we use the open
source package EZyRB~\cite{demo2018ezyrb}.

\section{Numerical results}
\label{sec:results}

Here we are going to present the results of the complete numerical
pipeline applied to the DTMB~5415 hull. Moreover we demonstrate the
improvements obtained using the proposed pipeline, called POD+AS,
with respect to the POD approach on the full parameter space. 

\begin{figure}[h!]
\centering
\includegraphics[trim=0 0 0 0, width=.8\textwidth]{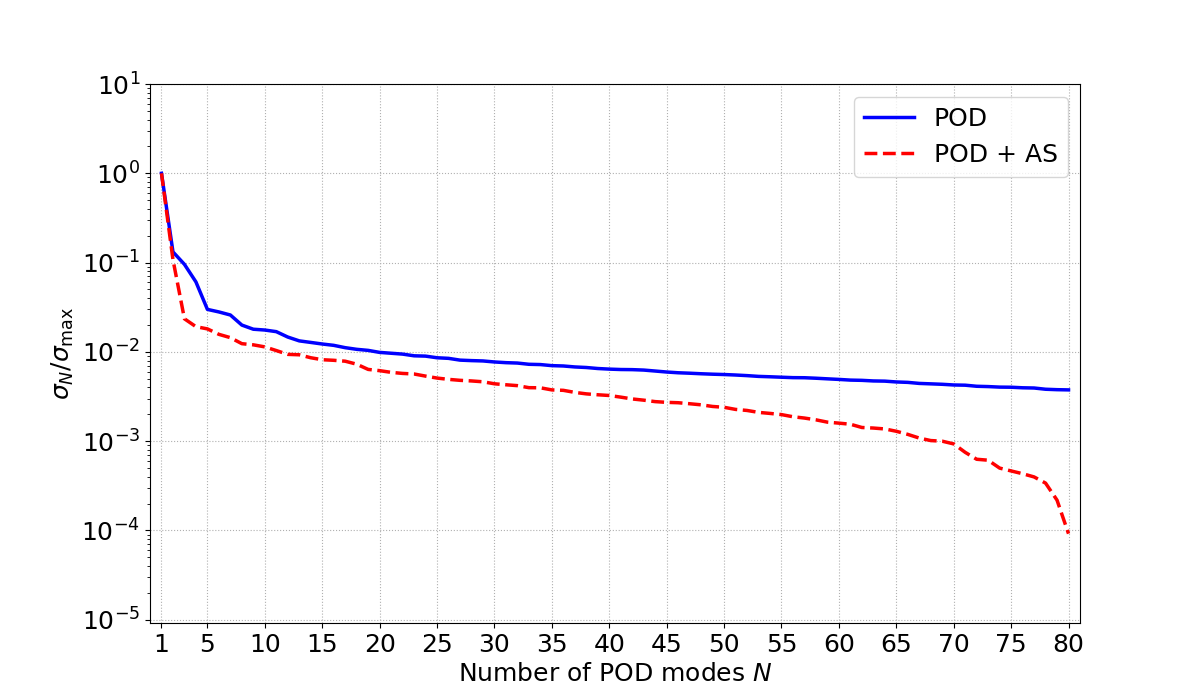}
\caption{POD singular values decay as a function of the number of
  modes. The blue line corresponds to the original sampling, while the
  red dotted line, called POD+AS approach, corresponds to the
  sampling along the active direction.}
\label{fig:svd}
\end{figure}

After generating $N_{\text{POD}} = 100$ deformed hulls, we perform the
high-fidelity simulations accelerated via the DMD algorithm. We
construct the snapshots matrix and compute the POD modes and the
corresponding eigenvalues for the construction of the reduced output
space. We compare this approach with the one proposed in this work
that exploits a preprocessing step with the finding of the active
subpace for the total drag resistance. With the $N_{\text{POD}}$
input/output couples, we individuate an eigenvector
$\mathbf{W}_1$ (compare Section~\ref{sec:active}) describing an active
subspace of dimension 1, and we sample the full space only along the
active direction described by this vector. After this second sampling
we collect a new set of high-fidelity simulations formed by
$N_{\text{POD+AS}} = 80$ snapshots. For this new snapshots matrix we
compute again the POD modes and eigenvalues, and we compare the two
approaches looking at the POD singular values decay. A faster decay
means a better approximation of the output fields for a fixed number
of modes. In Figure~\ref{fig:svd} the blue line shows the
singular values $\sigma_i$ divided by the first and greatest singular
value $\sigma_{\text{max}}$ for the sampling of the full parameter
space; with the dashed red line the POD singular values decay for the
POD+AS approach. A faster decay is observed, especially for the first
few modes.
This translates in an enhanced reduced order model, which exhibits a
better approximation of the solutions manifold, with respect to the
classical approach. 

Since we are relying on multidimensional interpolation to reconstruct
the solutions at untried parameters, having a new reduced parameter
space improves the creation of such interpolator. In the POD+AS
approach we have to interpolate a univariate function in
$N$-dimension, where $N$ is the number of POD modes we retain. In the
POD approach on the full parameter space instead, we have the same
number of modes to fit but a multivariate function depending on 5
input parameters, resulting in a difficult interpolation.

\section{Conclusions and perspectives}
\label{sec:the_end}

In this work we presented a nonintrusive numerical pipeline for shape
optimization of the bulbous bow of a benchmark hull. It comprises
automatic geometrical parametrization and morphing through FFD,
estimation of the total drag resistance via DMD using only few
snapshots of the time-dependent high fidelity simulations, the
reduction of the parameter space exploiting the AS property, and the
construction of a surrogate model with PODI for the real-time
evaluation of the many-query problem solved by an optimization
algorithm. We proved that the reduction of the parameter space can
further enhance the reduced order model creation. Moreover all this
parts of the pipeline can be used and integrated separately into an
existing computational workflow resulting in a great interest for
industrial applications.

\section*{Acknowledgements}
This work was partially performed in the context of the project SOPHYA -
``Seakeeping Of Planing Hull YAchts'' supported by Regione
FVG, POR-FESR 2014-2020, Piano Operativo Regionale Fondo Europeo per
lo Sviluppo Regionale, and partially supported by European Union Funding for
Research and Innovation --- Horizon 2020 Program --- in the framework
of European Research Council Executive Agency: H2020 ERC CoG 2015
AROMA-CFD project 681447 ``Advanced Reduced Order Methods with
Applications in Computational Fluid Dynamics'' P.I. Gianluigi
Rozza.

\end{document}